\documentclass[12pt,a4paper,twoside]{article}
\usepackage[english]{babel}
\usepackage{amssymb}
\usepackage{inputenc}
\newcommand{\ba}{\begin{array}}
\newcommand{\ea}{\end{array}}
\def\l {L^{0}}
\def\lt {L^{\infty}(\Omega, \mu)}
\usepackage{amssymb}
\begin{document}

\author{S. Albeverio $^{1},$ Sh. A. Ayupov $^{2},$ \ \ K. K.
Kudaybergenov  $^3$}
\title{\bf Derivations on   the Algebra of $\tau$-Compact Operators Affiliated with a   Type I von Neumann Algebra}

\maketitle

\begin{abstract}

Let $M$ be a type I von Neumann algebra with the center $Z,$ a
faithful normal semi-finite trace $\tau.$ Let $L(M, \tau)$ be the
algebra of all $\tau$-measurable operators affiliated with $M$ and
let $S_0(M, \tau)$ be the subalgebra in   $L(M, \tau)$ consisting
of all operators $x$ such that given any $\varepsilon>0$ there is
a projection  $p\in\mathcal{P}(M)$ with
$\tau(p^{\perp})<\infty,\,xp\in M$ and $\|xp\|<\varepsilon.$ We
prove that any $Z$-linear derivation of $S_0(M, \tau)$ is spatial
and generated by an element from  $L(M, \tau).$

\end{abstract}

\medskip
$^1$ Institut f\"{u}r Angewandte Mathematik, Universit\"{a}t Bonn,
Wegelerstr. 6, D-53115 Bonn (Germany); SFB 611, BiBoS; CERFIM
(Locarno); Acc. Arch. (USI),  e-mail: \emph{albeverio@uni-bonn.de}

$^2$ Institute of Mathematics, Uzbekistan Academy of Science, F.
Hodjaev str. 29, 700143, Tashkent (Uzbekistan), e-mail:
\emph{sh\_ayupov@mail.ru, e\_ayupov@hotmail.com,
mathinst@uzsci.net}

 $^{3}$ Institute of Mathematics, Uzbekistan
Academy of Science, F. Hodjaev str. 29, 700143, Tashkent
(Uzbekistan), e-mail: \emph{karim2006@mail.ru}



\newpage

\section*{\center 1. Introduction}

        The  present paper is devoted to  the description  of derivations
on a certain class of  algebras of $\tau$-measurable operators
affiliated with a type I von Neumann algebra.

It is known that all derivations on von Neumann algebras and more
general Banach algebras of operators are inner. Such kinds of
results were generalized for some classes of unbounded operator
algebras \cite{Sch} . In particular in \cite{AlAyup1}  it was
proved that any derivation on the non commutative Arens algebra
$L^{\omega}(M, \tau)$ associated with a von Neumann algebra $M$
with a faithful normal semi-finite trace  $\tau$ is spatial, and
if $\tau$ is finite then any derivation on   $L^{\omega}(M, \tau)$
is inner. Further if we consider the algebra  $L(M, \tau)$ of all
$\tau$-measurable operators affliliated with a type I von Neumann
algebra  $M,$ then a derivation $d$ on $L(M, \tau)$ is inner if
and only if it is $Z$-linear, where $Z$ is the center of $M$ (see
\cite{AlAyup}). Moreover there are examples of non $Z$-linear (and
hence non spatial and discontinuous in the measure topology)
derivations on $L(M, \tau)$ (see \cite{AlAyup}, \cite{Ber},
\cite{Kus}).

In this paper we study derivation on a subalgebra of $L(M,
\tau),$ namely on $$S_0(M, \tau)=\{x\in L(M, \tau):
\forall\varepsilon>0,\, \exists p,\,\tau(p^{\perp})<\infty,\,xp\in
M,\,\|xp\|<\varepsilon\},$$ which is an ideal in $L(M, \tau).$

In the particular case where $M=B(H)$ -- the algebra of all
bounded linear operators on a Hilbert space $H$ and $\tau=Tr$ --
the canonical trace,  $S_0(M, \tau)$ coincides with the ideal
$\mathcal{K}(H)$ of all compact operators on $H.$ In the general case the elements of $S_0(M, \tau)$ are called $\tau$-compact operators affiliated with the von Neumann algebra $M$ and the trace $\tau.$

Since $S_0(M, \tau)$ is an ideal in $L(M, \tau),$ any element
$a\in L(M, \tau)$ implements a derivation $d_a$ on $S_0(M, \tau)$
by
$$d_a(x)=ax-xa, \quad x\in S_0(M, \tau),$$
and such derivations are called spatial derivations. Moreover it is
clear that any spatial derivation $d_a$ is $Z$-linear, where $Z$
is the center of $M.$ The main result of the present paper states
the converse, i. e. that any $Z$-linear derivation on $S_0(M,
\tau)$ is spatial and implemented by an element of $L(M, \tau)$
for any type I von Neumann algebra $M.$ In particular if the
lattice of projections of $M$ is atomic, then any derivation on
$S_0(M, \tau)$ is automatically spatial.

In Section 2 we give some preliminaries from the theory of
lattice-normed modules and Kaplansky-Hilbert modules over the
algebra of measurable functions, and recall the main theorem from
\cite{AyupKud} which describes derivations on standard subalgebras
of the algebra of all $\l$-linear $\l$-bounded operators on a
Banach-Kantorovich space.

In Section 3 we give the main result which states that for any
type I von Neumann algebra $M$ with the center $Z,$ any $Z$-linear
derivation on the algebra $S_0(M, \tau)$  is spatial and
implemented by an element of $L(M, \tau).$ 

\begin{center} {\bf 2. Preliminaries}
\end{center}

 Let   $(\Omega,\Sigma,\mu)$  be a measurable space with
    a $\sigma$-finite measure $\mu,$ i. e. there is family $\{\Omega_{i}\}_{i\in
J}\subset\Sigma,\,0<\mu(\Omega_{i})<\infty,\,i\in J,$ such that
for any $A\in\Sigma,\,\mu(A)<\infty,$ there exists a countable subset
$J_{0 }\subset J$ and a set  $B$ with zero measure such that
$A=\bigcup\limits_{i\in J_{0}}(A\cap \Omega_{i})\cup B.$

 We
denote by  $\l=L^{0}(\Omega, \Sigma, \mu)$ the algebra of all
(classes of) complex measurable functions on $(\Omega, \Sigma,
\mu)$ equipped with the topology of convergence in measure. Then
$\l$ is a complete metrizable commutative regular algebra with the
unit $\textbf{1}$ given by $\textbf{1}(\omega)=1,
\,\omega\in\Omega.$

Denote by $\nabla$ the complete Boolean algebra of all idempotents
from $\l,$ i. e. $\nabla=\{\chi_{A}: A\in\Sigma\},$ where
$\chi_{A}$ is the characteristic function of the set  $A.$

A complex linear space  $E$ is said to be normed  by $\l$ if there
is
 a map  $\|\cdot\|:E\longrightarrow \l$ such that for any  $x,y\in E, \lambda\in
\mathbb{C},$ the following conditions are fulfilled:

 $\|x\|\geq 0; \|x\|=0\Longleftrightarrow x=0$; $\|\lambda x\|=|\lambda|\|x\|$;
  $\|x+y\|\leq\|x\|+\|y\|$.

 The  pair  $(E,\|\cdot\|)$ is called a lattice-normed  space over $\l.$
A lattice-normed  space  $E$ is called  $d$-decomposable, if for
any $x\in E$ with  $\|x\|=\lambda_{1}+\lambda_{2},$ where $
\lambda_{1}, \lambda_{2} \in\l,\,\lambda_{1}\lambda_{2}=0,$ there
exists $x_{1}, x_{2}\in E$ such that $x=x_{1}+x_{2}$ and
$\|x_{i}\|=\lambda_{i},\,i=1,2$. A net $(x_{\alpha})$ in $E$ is
$(bo)$-converging to $x\in E$, if $\|x_{\alpha}-x\|\rightarrow0$
$\mu$-almost everywhere in
 $\l.$
 A lattice-normed  space $E$ which is $d$-decomposable and  complete with respect to the  $(bo)$-convergence
  is called a
 \emph{Banach-Kantorovich space}.

 It is known  that every Banach-Kantorovich space $E$ over
 $\l$ is a module over $\l$ and  $\|\lambda x\|=|\lambda|\|x\|$
 for all $\lambda\in\l,\, x\in E$ (see \cite{Kusr}).

Any Banach-Kantorovich space  $E$ over $\l$ is orthocomplete, i.
e. given any net  $(x_{\alpha})\subset E$ and a partition of the
unit $(\pi_{\alpha})$ in $\nabla$ the series
$\sum\limits_{\alpha}\pi_{\alpha}x_{\alpha}$ $(bo)$-converges in
$E.$

 A module  $F$ over  $\l$ is said to be  finite-generated, if there are
       $x_{1},x_{2},...,x_{n}$ in $F$           for any  $x\in F$ there exists  $\lambda_{i}\in
       \l\, (i=\overline{1,n})$ such that $x=\lambda_{1}x_{1}+...+\lambda_{n}x_{n}.$
       The elements $x_{1},x_{2},...,x_{n}$ are called  generators of $F.$
       We denote by $d(F)$ the minimal number of  generators of $F.$
       A  module  $F$ over  $\l$ is  called $\sigma$-finite-generated,
       if there exists
       a partition $(\pi_{\alpha})_{\alpha\in A}$ of the unit in $\nabla$ such that
       $\pi_{\alpha}F$ is finite-generated for any $\alpha.$
         A finite-generated module  $F$ over  $\l$ is called homogeneous of type  $n$, if for every
          nonzero $e\in\nabla$ we have  $n=d(eF)$.

 Let     $E$ be a Banach-Kantorovich space  over $\l.$
              If $(u_{\alpha})_{\alpha\in A}\subset E$ and
$(\pi_{\alpha})_{\alpha\in A}$
              is a partition of the unit in $\nabla$, then the
              series
      $\sum\limits_{\alpha}\pi_{\alpha}u_{\alpha}$ $(bo)$-converges in
      $E$ and  its sum is called the mixing of
       $(u_{\alpha})_{\alpha\in A}$ with respect to  $(\pi_{\alpha})_{\alpha\in A}.$ We denote this sum
             by $\textmd {mix}(\pi_{\alpha}u_{\alpha})$.
      A subset $K\subset E$ is called cyclic, if
       $\textmd {mix}(\pi_{\alpha}u_{\alpha})\in K$ for each $(u_{\alpha})_{\alpha\in A}\subset K$ and
       any partition  of the unit $(\pi_{\alpha})_{\alpha\in A}$ in $\nabla.$
       For every directed set $A$  denote by  $\nabla (A)$
       the set of all  partitions of the unit in $\nabla,$ which are  indexed by elements of the set
       $A.$  More precisely,
       $$
       \nabla(A)=\{\nu : A\rightarrow \nabla : (\forall\alpha,\beta
       \in A)(\alpha\neq\beta\rightarrow\nu(\alpha)\wedge\nu(\beta)=0)\wedge\bigvee_{\alpha\in
       A}\nu(\alpha)=\textbf{1}\}.
       $$
For $\nu_{1}, \nu_{2} \in\nabla(A)$ we put
$\nu_{1}\leq\nu_{2}\leftrightarrow \forall\alpha,\beta\in A,\,
      (\nu_{1}(\alpha)\wedge\nu_{2}(\beta)\neq0\rightarrow\alpha\leq\beta).$
           Then $\nabla (A)$ is a directed set.
       Let  $(u_{\alpha})_{\alpha\in A}$ be a net in $E.$ For
       every  $\nu\in \nabla (A)$  we put  $u_{\nu}=\textmd {mix}(\nu(\alpha)u_{\alpha})$
       and    obtain a new net $(u_{\nu})_{\nu\in \nabla (A)}.$
        Every  subnet of the net
              $(u_{\nu})_{\nu\in \nabla (A)} $  is called a cyclic subnet of
            the original net   $(u_{\alpha})_{\alpha\in A}.$

      \textbf{  Definition} \cite {Kusr}.  A subset $K\subset E$
       is called  \emph{cyclically
       compact,} if
        $K$ is cyclic and   every net in $K$ has a  cyclic subnet that $(bo)$-converges
         to some point of  $K.$         A subset in $E$ is called \emph{relatively cyclically
       compact} if it is contained in a cyclically
       compact set.

Let  $\mathcal{K}$ be a module over $\l$. A map $\langle
\cdot,\cdot\rangle:\mathcal{K}\times \mathcal{K}\rightarrow\l$ is
called  an $\l$-valued inner product, if for all $x,y,z\in
\mathcal{K},\,\lambda\in\l,$ the following conditions are
fulfilled: $\langle x,x\rangle\geq0$; $\langle
x,x\rangle=0\Leftrightarrow x=0$; $\langle
x,y\rangle=\overline{\langle y,x\rangle}$; $\langle \lambda
x,y\rangle=\lambda\langle x,y\rangle$; $\langle
x+y,z\rangle=\langle x,z\rangle+\langle y,z\rangle$.

If $\langle \cdot,\cdot\rangle:\mathcal{K}\times
\mathcal{K}\rightarrow\l$ is an $\l$-valued inner product, then $
\|x\|=\sqrt{\langle x,x \rangle} $ defines  an $\l$-valued norm on
$\mathcal{K}.$ The  pair $(\mathcal{K},\langle
\cdot,\cdot\rangle)$ is called a \emph{Kaplansky-Hilbert module}
over $\l,$ if $(\mathcal{K},\|\cdot\|)$ is a Banach-Kantorovich
space over  $\l$ (see \cite{Kusr}).

Let  $X$ be a Banach space.  A map  $s:\Omega\rightarrow X$ is
called a simple, if
$s(\omega)=\sum\limits_{k=1}^{n}\chi_{A_{k}}(\omega)c_k,$ where
$A_k\in\Sigma, A_i\cap A_j=\emptyset, \,i\neq j,\, \,c_k\in
X,\,k=\overline{1, n},\, n\in\mathbb{N}.$ A map
$u:\Omega\rightarrow X$ is said to be measurable, if there is a
sequence  $(s_n)$ of simple maps such that
$\|s_n(\omega)-u(\omega)\|\rightarrow0$  almost everywhere on any
$A\in\sum$ with $\mu(A)<\infty.$

Let $\mathcal{L}(\Omega, X)$ be the set of all measurable maps
from  $\Omega$ into $X,$ and let $L^{0}(\Omega, X)$ denote the
factor\-ization of this set with respect to equality  almost
everywhere. Denote by $\hat{u}$ the equivalence class from
$L^{0}(\Omega, X)$ which contains the measurable map $u\in
\mathcal{L}(\Omega, X).$ Further we shall identity the element
$u\in \mathcal{L}(\Omega, X)$ with the class $\hat{u}.$ Note that
the function  $\omega \rightarrow \|u(\omega)\|$
       is measurable for any $u\in \mathcal{L}(\Omega, X).$ The equivalence class containing the function
              $\|u(\omega)\|$ is denoted by
       $\|\hat{u}\|$. For  $\hat{u}, \hat{v}\in L^{0}(\Omega, X), \lambda\in\l$ put
$\hat{u}+\hat{v}=\widehat{u(\omega)+v(\omega)},
\lambda\hat{u}=\widehat{\lambda(\omega) u(\omega)}.$

It is known  \cite{Kusr} that $(L^{0}(\Omega, X), \|\cdot\|)$ is
a Banach-Kantorovich space over $\l.$

Put $L^{\infty}(\Omega, X)=\{x\in L^{0}(\Omega, X):\|x\|\in
L^{\infty}(\Omega)\}.$  Then $L^{\infty}(\Omega, X)$ is a Banach
space with respect the norm
$\|x\|_{\infty}=\|\|x\|\|_{L^{\infty}(\Omega)}.$

If $H$ is a Hilbert space, then $L^{0}(\Omega, H)$ can be equipped
with an  $\l$-valued inner product $\langle x,
y\rangle=\widehat{(x(\omega), y(\omega))},$ where $(\cdot, \cdot)$
is the inner product on $H.$

Then $(L^{0}(\Omega, H),\langle \cdot,\cdot\rangle)$ is a
Kaplansky-Hilbert module over $\l.$

Let  $E$ be a  Banach-Kantorovich space over $\l.$
 An operator   $T: E\rightarrow E$
 is called  $\l$-linear if  $T(\lambda_1 x_1 +\lambda_2 x_2)=\lambda_1 T(x_1)+\lambda_2
 T(x_2)$
for all  $\lambda_1, \lambda_2\in \l ,x_1, x_2\in E.$ An
$\l$-linear  operator
 $T:E\rightarrow E$ is called $\l$-bounded if there exists an element $c\in\l$
 such that $\|T(x)\|\leq
 c\|x\|$ for any $x\in E.$  For an  $\l$-bounded linear
operator $T$  we put $\|T\|=\sup\{\|T(x)\|:\|x\|\leq \textbf{1}
\}.$

An $\l$-linear operator $T:E\rightarrow E$ is called
finite-generated ($\sigma$-finite-generated, homogeneous of type
$n$) if $T(E)=\{T(x):x\in E\}$ is a finite-generated (respectively
$\sigma$-finite-generated, homogeneous of type $n$) submodule in
$E.$

 An $\l$-linear
operator $T:E\rightarrow E$
 is called cyclically compact, if for every bounded set $B$ in $E$ the set
 $T(B)$ is relatively cyclically compact  in $E.$

We denote by  $B(E)$   the  algebra of all
 $\l$-linear $\l$-boun\-ded operators
on  $E$
 and  $\mathcal{F}(E)$  be the set of all
finite-generated $\l$-linear $\l$-bounded operators on $E.$

An algebra  $\mathcal{U}\subset B(E)$ is called \emph{standard}
over  $\l,$ if $\mathcal{U}$ is a submodule in $B(E)$ and
$\mathcal{F}(E)\subset\mathcal{U}.$

The examples of standard subalgebras are given by $B(E),$
$\mathcal{F}(E)$ and the space of all cyclically compact operators
on $E.$

 \textbf{Theorem 2.1} \cite {AyupKud}. \emph{Let $\mathcal{U}$ be a  standard algebra in    $B(E)$
and let $\delta:\mathcal{U}\rightarrow B(E)$ be an $\l$-linear
derivation.
 Then there is  $T\in B(E)$ such that $\delta(A)=TA-AT$ for all} $A\in \mathcal{U}.$

 \begin{center} {\bf 3. The main result}
\end{center}

Let  $B(H)$ be the algebra of all bounded linear operators on a
Hilbert space  $H$ and let  $M$ be a von Neumann algebra in $B(H)$
with a faithful normal semi-finite trace $\tau.$ Denote by
$\mathcal{P}(M)$ the lattice of projections in  $M.$

A linear subspace  $D$ in  $H$ is said to be affiliated with  $M$
(denotes as  $D\eta M$), if $u(D)\subset D$ for any unitary
operator  $u$ from the commutant $$M'=\{y'\in B(H):xy'=y'x,
\,\forall x\in M\}$$ of the algebra $M.$

A linear operator  $x$ on  $H$ with domain  $D(x)$ is said to be
affiliated with  $M$ (denoted as  $x\eta M$) if  $u(D(x))\subset
D(x)$ and $ux(\xi)=xu(\xi)$ for all $u\in M',$ $\xi\in D(x).$

A linear subspace  $D$ in $H$ is called  $\tau$-dense, if

1) $D\eta M;$

2) given any  $\varepsilon>0$ there exists a projection
$p\in\mathcal{P}(M)$ such that  $p(H)\subset D$ and
$\tau(p^{\perp})\leq\varepsilon.$

A closed linear operator  $x$ is said to be $\tau$-measurable (or
totally measurable) with respect to the von Neumann algebra  $M,$
if  $x\eta M$ and $D(x)$ is  $\tau$-dense in $H.$

 We will denote by $L(M,
\tau)$ the set of all $\tau$-measurable operators affiliated with
$M.$ Let $\|\cdot\|_{M}$ stand for the uniform norm in $M.$ The
\emph{measure topology,} $t_{\tau},$ in $L(M, \tau)$ is the one
given by the following system of neighborhoods of zero:
$$V(\varepsilon, \delta)=\{x\in L(M, \tau): \exists e\in\mathcal{P}(M), \tau(e^{\perp})\leq\delta, xe\in
M,  \|xe\|_{M}\leq\varepsilon\},$$ where $\varepsilon>0,
\delta>0.$

It is known \cite{Nel} that 
$L(M, \tau)$ equipped with the measure topology is a complete metrizable topological $\ast$-algebra.

In the algebra   $L(M, \tau)$ consider the subset  $S_0(M, \tau)$
of all operators $x$ such that given any $\varepsilon>0$ there is
a projection  $p\in\mathcal{P}(M)$ with
$\tau(p^{\perp})<\infty,\,xp\in M$ and $\|xp\|<\varepsilon.$  Following \cite{Str} let us call the elements of $S_0(M, \tau)$ \emph{$\tau$-compact operators} affiliated with $M.$ It is
known \cite{Yea}, \cite{ChilLit} that  $S_0(M, \tau)$ is a
$\ast$-subalgebra in  $L(M, \tau)$ and an $M$-bimodule, i. e.  $ax, xa\in
S_0(M, \tau)$ for all $x\in S_0(M, \tau)$ and $a\in M.$ It is
clear that if the trace  $\tau$ is finite then $S_0(M, \tau)=L(M,
\tau).$

The following properties of the algebra $S_0(M, \tau)$ of $\tau$-compact operators are known \cite{Str}, \cite{Bik}, but the proof is included for sake of completeness.

\textbf{Proposition 3.1.} \emph{Let  $M$ be a von Neumann algebra
with a faithful normal semi-finite trace $\tau.$ Then}

1) $L(M, \tau)=M+S_0(M, \tau);$

2) $S_0(M, \tau)$ \emph{is an ideal in} $L(M, \tau).$

Proof. Let  $x\in L(M, \tau).$ Take a projection  $p\in M$ such
that  $\tau(p^{\perp})<\infty$ and $xp\in M.$ Put $x_1=xp$ and
$x_2=xp^{\perp}.$ Since $\tau(p^{\perp})<\infty$ and $x_2p=0,$
then $x_2\in S_0(M, \tau).$ Therefore, any element from $L(M,
\tau)$ can be a represented as $x=x_1+x_2,$ where $x_1\in M,$
$x_2\in S_0(M, \tau).$ Since $S_0(M, \tau)$ is a module over $M,$
then from the equality  $L(M, \tau)=M+S_0(M, \tau)$ it follows
that $S_0(M, \tau)$ is an ideal in $L(M, \tau).$ The proof is
complete. $\blacksquare$

Since  $S_0(M, \tau)$ is an ideal in $L(M, \tau),$ any
element  $a\in L(M, \tau)$ implements a derivation on the algebra
$S_0(M, \tau)$ by the formula
$$d(x)=ax-xa,\quad x\in S_0(M, \tau),$$
which is $Z$-linear, $Z$ being the center of $M.$

 The main aim of the present work is to prove the converse, i. e. any $Z$-linear derivation
 on  $S_0(M,\tau)$ is spatial and implemented by an element of
$L(M, \tau).$

 Let
$\lt\bar{\otimes}B(H)$ be the tensor product of von  Neumann
algebra  $\lt$ and $B(H),$ with the trace $\tau=\mu\otimes Tr,$
where $Tr$ is the canonical trace for operators in $B(H)$  (with
its natural domain).

Denote by $L^{0}(\Omega, B(H))$ the space of equivalence classes
of measurable maps from  $\Omega$ into $B(H).$ Given  $\hat{u},
\hat{v}\in L^{0}(\Omega, B(H))$ put
$\hat{u}\hat{v}=\widehat{u(\omega) v(\omega)},
\hat{u}^{\ast}=\widehat{u(\omega)^{\ast}}.$

Define $$L^{\infty}(\Omega, B(H))=\{x\in L^{0}(\Omega,
B(H)):\|x\|\in L^{\infty}(\Omega)\}.$$ The space
$(L^{\infty}(\Omega, B(H)), \|\cdot\|_{\infty})$ is a Banach
*-algebra.

 It is known  \cite{Tak} that the algebra  $\lt\bar{\otimes}B(H)$
is *-isomorphic to the algebra  $L^{\infty}(\Omega, B(H)).$

Note also that
$$\tau(x)=\int\limits_{\Omega}Tr(x(\omega))\,d\mu(\omega).$$

Further we shall identity the algebra  $\lt\bar{\otimes}B(H)$ with
the algebra $L^{\infty}(\Omega, B(H)).$

Denote by $B(L^{0}(\Omega, H))$ (resp. $B(L^{\infty}(\Omega, H))$)
the algebra of all  $\l$-linear and $\l$-bounded  (resp.
$L^{\infty}(\Omega)$-linear and $L^{\infty}(\Omega)$-bounded)
operators on $L^{0}(\Omega, H)$ (resp. $L^{\infty}(\Omega, H)$).

Given any $f\in L^{\infty}(\Omega, B(H))$  consider the element
$\Psi(f)$ from $B(L^{\infty}(\Omega, H))$ defined by
 $$\Psi(f)(x)=\widehat{f(\omega)(x(\omega))},\quad x\in L^{\infty}(\Omega, H).$$

Then the correspondence  $f\rightarrow \Psi(f)$ gives an isometric
*-isomorphism between the algebras  $L^{\infty}(\Omega, B(H))$ and
$B(L^{\infty}(\Omega, H))$  (see \cite{Kusr}).

Since $L^{\infty}(\Omega, B(H))$ is $(bo)$-dense in $L^{0}(\Omega,
B(H))$ and  $B(L^{\infty}(\Omega, H))$ is $(bo)$-dense in
$B(L^{0}(\Omega, H)),$ the  *-isomorphism  $\Psi$ can be uniquely
extended to a *-isomorphism between  $L^{0}(\Omega, B(H))$ and
$B(L^{0}(\Omega, H)).$

It is known \cite{AlAyup}, that the algebra
$L(\lt\bar{\otimes}B(H), \tau)$ of all $\tau$-measurable operators
affiliated with the von Neumann algebra $\lt\bar{\otimes}B(H)$ is
$\l$-linear *-isomorphic with the algebra $B(L^{0}(\Omega, H)).$

Therefore one has the following relations for the algebras
mentioned above:
$$
\begin{array}{ccccc}
 L^{\infty}(\Omega)\bar{\otimes} B(H)  &  \cong& L^{\infty}(\Omega, B(H)) & \cong & B(L^{\infty}(\Omega,H)) \\
  \cap &  & \cap &  & \cap \\
  L(L^{\infty}(\Omega)\bar{\otimes} B(H)),\tau) & \cong & L^{0}(\Omega, B(H))  & \cong & B(L^{0}(\Omega,H)). \\
\end{array}
$$

\textbf{Proposition 3.2}. \emph{Let
$p\in\mathcal{P}(L^{\infty}(\Omega, B(H)))$ and $\tau(p)<\infty.$
Then  $p$ is $\sigma$-finite-generated and in particular is
cyclically compact. }

 Proof. Since $\tau(p)=\int\limits_{\Omega}Tr(p(\omega))\,d\mu(\omega)$
we have that  $Tr(p(\omega))<\infty$ for almost all
$\omega\in\Omega.$ In the algebra  $B(H)$ any projection with
finite trace is finite dimensional, thus  $p(\omega)$ is a finite
dimensional projection for almost all  $\omega\in\Omega.$ By
(\cite{GK}, Theorem 2)  $p$ is $\sigma$-finite-generated and thus
 $p$ is cyclically compact. The proof is complete. $\blacksquare$

\textbf{Proposition  3.3}. \emph{If  $x\in
S_0(L^{\infty}(\Omega)\bar{\otimes} B(H), \tau)$ then  $x$ is
cyclically compact.}

Proof. If  $x\in S_0(L^{\infty}(\Omega)\bar{\otimes} B(H), \tau)$
then given any  $\varepsilon>0$ there is a projection
$p_{\varepsilon}\in\mathcal{P}(L^{\infty}(\Omega)\bar{\otimes}
B(H))$ such that
$$\tau(p^{\perp}_{\varepsilon})<\infty,\,xp_{\varepsilon}\in L^{\infty}(\Omega)\bar{\otimes} B(H),\,
\|xp_{\varepsilon}\|_{\infty}<\varepsilon.$$

By Proposition 3.2 $p^{\perp}_{\varepsilon}$ is cyclically compact
and therefore  $xp^{\perp}_{\varepsilon}$ is also cyclically
compact. From $\|xp_{\varepsilon}\|_{\infty}<\varepsilon$ it
follows that $\|xp_{\varepsilon}\|\leq\varepsilon\textbf{1}.$
Therefore
$\|x-xp^{\perp}_{\varepsilon}\|\leq\varepsilon\textbf{1},$ i. e.
$x$ is the  $(bo)$-limit of cyclically compact operators and thus
 $x$ is also cyclically compact. The proof is complete. $\blacksquare$

The converse assertion for Proposition 3.3 is not true in general.
Indeed, let  $\mu(\Omega)=+\infty$ and $\dim H<\infty.$ Then
$L(L^{\infty}(\Omega)\bar{\otimes} B(H), \tau)$ is *-isomorphic to
the algebra  of $n\times n$ matrices over  $\l.$ Therefore any operator
from  $L(L^{\infty}(\Omega)\bar{\otimes} B(H), \tau)$ is
cyclically compact because it acts on the finite-generated module
over  $\l.$ In particular the identity $e$ in
$L(L^{\infty}(\Omega)\bar{\otimes} B(H), \tau)$ is cyclically
compact. But   $e\notin S_0(L^{\infty}(\Omega)\bar{\otimes} B(H),
\tau)$ because $\mu(\Omega)=\infty.$

Let   $\mbox{mix}(S_0(L^{\infty}(\Omega)\bar{\otimes} B(H),
\tau))$ be the cyclic hull of the set
$S_0(L^{\infty}(\Omega)\bar{\otimes} B(H), \tau),$ i. e. it 
consists of all elements of the form
$x=(bo)-\sum\limits_{\alpha}\pi_\alpha x_\alpha, $ where
$(\pi_\alpha)$ is a partition of the unit in  $\nabla,$
$(x_\alpha)\subset S_0(L^{\infty}(\Omega)\bar{\otimes} B(H),
\tau).$

Since  $S_0(L^{\infty}(\Omega)\bar{\otimes} B(H), \tau)$ is a
module over  $L^{\infty}(\Omega)$ and
$\l=\mbox{mix}(L^{\infty}(\Omega)),$ we have that
$\mbox{mix}(S_0(L^{\infty}(\Omega)\bar{\otimes} B(H), \tau))$ is a
module over $\l.$

\textbf{Proposition 3.4}.
\emph{$\mbox{mix}(S_0(L^{\infty}(\Omega)\bar{\otimes} B(H),
\tau))$ is a standard algebra in
$L(L^{\infty}(\Omega)\bar{\otimes} B(H), \tau).$}

Proof. First suppose that the measure  $\mu$ is finite. Consider a
finite-generated operator  $x$ from the algebra
$L(L^{\infty}(\Omega)\bar{\otimes} B(H), \tau).$ Let $p$ be the
orthogonal projection onto the image of $x$ and let $n$  -- be the
number of its generators. By (\cite{GK}, Theorem 2)
$Tr(p(\omega))=\dim p(\omega)\leq n$  for almost all
$\omega\in\Omega.$ Therefore $\tau(p)=\int\limits_{\Omega}
Tr(p(\omega))d\mu(\omega)\leq n\mu(\Omega),$ i. e.
$\tau(p)<\infty.$

It is clear that  $xp^{\perp}=0.$ Thus $\tau(p)<\infty$ and
$xp^{\perp}=0,$ i. e. $x\in S_0(L^{\infty}(\Omega)\bar{\otimes}
B(H), \tau).$

Now suppose that  $\mu$ is  $\sigma$-finite and  $x$ is a
finite-generated operator from $L(L^{\infty}(\Omega)\bar{\otimes}
B(H), \tau).$ Since the measure  $\mu$  is $\sigma$-finite, there
exists a partition of the unit  $(e_\alpha)$ in $\nabla$ such that
$e_\alpha=\chi_{A_\alpha}, A_\alpha\in\Sigma,
\mu(A_\alpha)<\infty.$ From the above it follows $e_\alpha x\in
e_\alpha S_0(L^{\infty}(\Omega)\bar{\otimes} B(H), \tau)$ and
therefore $x=(bo)-\sum\limits_{\alpha}e_\alpha x$ belongs to
$\mbox{mix}(S_0(L^{\infty}(\Omega)\bar{\otimes} B(H), \tau)).$
Thus  $\mbox{mix}(S_0(L^{\infty}(\Omega)\bar{\otimes} B(H),
\tau))$ is a standard algebra. The proof is complete.
$\blacksquare$

\textbf{Proposition 3.5}. \emph{Any $L^{\infty}(\Omega)$-linear
derivation  $d$ of the algebra
$S_0(L^{\infty}(\Omega)\bar{\otimes} B(H), \tau)$ is spatial and}
$$d(x)=ax-xa, \quad x\in S_0(L^{\infty}(\Omega)\bar{\otimes} B(H), \tau),\eqno (1)$$
\emph{for an appropriate} $a\in L(L^{\infty}(\Omega, B(H)),
\tau).$

Proof. Let $d$ be a $L^{\infty}(\Omega)$-linear derivation of the
algebra $S_0(L^{\infty}(\Omega)\bar{\otimes} B(H), \tau).$ Let us
show that  $d$ can be extended onto the algebra
$\mbox{mix}(S_0(L^{\infty}(\Omega)\bar{\otimes} B(H), \tau)).$ By
definition any element of
$\mbox{mix}(S_0(L^{\infty}(\Omega)\bar{\otimes} B(H), \tau))$ has
the form  $$x=(bo)-\sum\limits_{\alpha}\pi_\alpha x_\alpha,
$$ where  $(\pi_\alpha)$ is a partition of the unit in $\nabla,$
$(x_\alpha)\subset S_0(L^{\infty}(\Omega)\bar{\otimes} B(H),
\tau).$

Put
$$\tilde{d}(x)=(bo)-\sum\limits_{\alpha}\pi_\alpha
d(x_\alpha).$$

Straightforward arguments show that  $\tilde{d}$ is a well-defined
derivation on the algebra
$\mbox{mix}(S_0(L^{\infty}(\Omega)\bar{\otimes} B(H), \tau)).$

Let us prove that  $\tilde{d}$ is $\l$-linear. Let $\lambda\in \l$
and $x\in\mbox{mix}(S_0(L^{\infty}(\Omega)\bar{\otimes} B(H),
\tau)).$ Take a partition of the unit $(e_{\alpha})$ in $\nabla$
such that $e_{\alpha}\lambda\in L^{\infty}(\Omega),$ $e_\alpha
x\in S_0(L^{\infty}(\Omega)\bar{\otimes} B(H), \tau)$ for all
$\alpha.$ Since  $d$ is  $L^{\infty}(\Omega)$-linear we have
$d(e_{\alpha}\lambda x)=d(e_{\alpha}\lambda
e_{\alpha}x)=e_{\alpha}\lambda d(e_{\alpha} x).$ Therefore
$\tilde{d}(\lambda x)=(bo)-\sum\limits_{\alpha}e_\alpha
d(e_\alpha\lambda x)=(bo)-\sum\limits_{\alpha}e_\alpha \lambda
d(e_{\alpha} x)=\lambda \tilde{d}(x),$ i. e.  $\tilde{d}(\lambda
x)=\lambda \tilde{d}(x).$

Since  $\mbox{mix}(S_0(L^{\infty}(\Omega)\bar{\otimes} B(H),
\tau))$ is a standard algebra Theorem 2.1 implies that the
derivation  $\tilde{d}$ and, in particular, the derivation $d$ is
of the form  (1). The proof is complete. $\blacksquare$

Recall that a von Neumann algebra $M$ is an algebra of \emph{type
I} if it is isomorphic to a von Neumann algebra with an abelian
commutant.

It is well-known \cite{Tak} that if   $M$ is a type I von Neumann
algebra
 then there is a unique (cardinal-indexed) orthogonal
family of projections $(q_{\alpha})_{\alpha\in
I}\subset\mathcal{P}(M)$ with $\sum\limits_{\alpha\in
I}q_{\alpha}=\textbf{1}$ such that  $q_{\alpha}M$ is isomorphic to
the tensor product of an abelian von Neumann algebra
$L^{\infty}(\Omega_{\alpha}, \mu_{\alpha})$ and $B(H_{\alpha})$
with $\dim H_{\alpha}=\alpha,$ i. e.
$$M\cong\sum\limits_{\alpha}^{\oplus}L^{\infty}(\Omega_{\alpha}, \mu_{\alpha})\bar{\otimes}B(H_{\alpha}).$$

Consider the faithful normal semi-finite trace $\tau$ on $M,$
defined as
$$\tau(x)=\sum\limits_{\alpha}\tau_{\alpha}(x_{\alpha}),\quad x=(x_{\alpha})\in M,\,x\geq0,$$
where $\tau_{\alpha}=\mu_{\alpha}\otimes Tr_{\alpha}.$

Now we can prove the main result of the present paper.

\textbf{Theorem 3.6}. \emph{If  $M$ is a von Neumann algebra of
type I with the center $Z,$ then any $Z$-linear derivation on the
algebra $S_0(M, \tau)$ is spatial and  implemented by an element
of $L(M, \tau).$}

In order to prove the theorem we need several auxiliary results.

Let

 $$\prod\limits_{\alpha}L(L^{\infty}(\Omega_{\alpha},
\mu_{\alpha})\bar{\otimes}B(H_{\alpha}), \tau_{\alpha})$$ be the
topological (Tychonoff) product of the spaces
$L(L^{\infty}(\Omega_{\alpha},
\mu_{\alpha})\bar{\otimes}B(H_{\alpha}), \tau_{\alpha}).$

  Then (see \cite{ChilLit}) we have the topological embedding
$$L(M, \tau)\subset\prod\limits_{\alpha}L(L^{\infty}(\Omega_{\alpha},
\mu_{\alpha})\bar{\otimes}B(H_{\alpha}), \tau_{\alpha}).$$

Denote by $Z_0$ the center of the algebra  $L(M, \tau).$ Then
$Z_0$ is *-isomorphic with the algebra of all $\tau'$-measurable
operators affiliated with the abelian von Neumann algebra
$\sum\limits_{\alpha}^{\oplus}L^{\infty}(\Omega_{\alpha},
\mu_{\alpha}),$ where the trace  $\tau'$ is defined by
$$\tau'(f)=\sum\limits_{\alpha}\int\limits_{\Omega_{\alpha}}f_{\alpha}d\mu_{\alpha},\,f\geq0.$$

Let $\Phi_{\alpha}$ be a  *-isomorphism between the algebras
$L(L^{\infty}(\Omega_{\alpha},
\mu_{\alpha})\bar{\otimes}B(H_{\alpha}), \tau_{\alpha})$ and
$B(L^{0}(\Omega_{\alpha}, H_{\alpha})).$ Given
$x=(x_{\alpha})\in\prod\limits_{\alpha}
L(L^{\infty}(\Omega_{\alpha},
\mu_{\alpha})\bar{\otimes}B(H_{\alpha}), \tau_{\alpha})$ put
$$\|x\|=(\|\Phi_{\alpha}(x_{\alpha})\|_{\alpha}),$$
where
  $\|\cdot\|_{\alpha}$  is  the norm on
$B(L^{0}(\Omega_{\alpha}, H_{\alpha})).$

Then an element   $x=(x_{\alpha})\in\prod\limits_{\alpha}
L(L^{\infty}(\Omega_{\alpha},
\mu_{\alpha})\bar{\otimes}B(H_{\alpha}), \tau_{\alpha})$ belongs
to $L(M, \tau)$ if and only if $\|x\|\in Z_{0}$ (see
\cite{AlAyup}).

 We consider on
$\prod\limits_{\alpha}L^{0}(\Omega_\alpha, H_\alpha)$ the
$\prod\limits_{\alpha}L^{0}(\Omega_\alpha)$-valued norm defined
by
$$\|\varphi\|=(\|\varphi_\alpha\|_{L(\Omega_\alpha, H_\alpha)}).$$

Then $\prod\limits_{\alpha}L^{0}(\Omega_\alpha, H_\alpha)$ is a
Banach-Kantorovich space over
$\prod\limits_{\alpha}L^{0}(\Omega_\alpha).$

Set
$$\oplus_{\alpha}L(\Omega_\alpha, H_\alpha)\equiv\{(\varphi_\alpha)\in\prod\limits_{\alpha}L^{0}(\Omega_\alpha, H_\alpha):
(\|\varphi_\alpha\|_{L(\Omega_\alpha, H_\alpha)})\in Z_0\}.$$
Since  $Z_0$ is a solid subalgebra in
$\prod\limits_{\alpha}L^{0}(\Omega_\alpha)$ it follows that
$(\oplus_{\alpha}L^{0}(\Omega_\alpha, H_\alpha), \|\cdot\|)$ is a
Banach-Kantorovich space over $Z_0.$

Let  $B(\oplus_{\alpha}L^{0}(\Omega_\alpha, H_\alpha))$ be  the
algebra of all $Z_0$-linear $Z_0$-bounded operators on
$\oplus_{\alpha}L^{0}(\Omega_\alpha, H_\alpha).$

Set
$$\oplus_{\alpha}B(L(\Omega_\alpha, H_\alpha))\equiv\{(x_\alpha)\in\prod\limits_{\alpha}B(L^{0}(\Omega_\alpha, H_\alpha)):
(\|x_\alpha\|_{B(L(\Omega_\alpha, H_\alpha))})\in Z_0\}.$$

It is clear that  $B(\oplus_{\alpha}L^{0}(\Omega_\alpha,
H_\alpha))$ is *-isomorphic to
$\oplus_{\alpha}B(L^{0}(\Omega_\alpha, H_\alpha)).$

\textbf{Lemma 3.7}. \emph{The algebra   $L(M, \tau)$ is
*-isomorphic to $B(\oplus_{\alpha}L^{0}(\Omega_\alpha,
H_\alpha)).$}

Proof. Let $\Phi_\alpha$ be a  *-isomorphism between
$L(L^{\infty}(\Omega_{\alpha},
\mu_{\alpha})\bar{\otimes}B(H_{\alpha}), \tau_\alpha)$ and
$B(L^{0}(\Omega_\alpha, H_\alpha)).$

Put $$\Phi(x)=(\Phi_\alpha(x_\alpha)),\quad  x\in L(M, \tau).$$
Since $x\in L(M, \tau)$ exactly means that  $\|x\|\in Z_0,$ and
$x'\in \oplus_{\alpha}B(L^{0}(\Omega_\alpha, H_\alpha))$ means
$\|x'\|\in Z_0,$  these imply that $\Phi$ is a *-isomophism
between  $L(M, \tau)$ and $ \oplus_{\alpha}B(L^{0}(\Omega_\alpha,
H_\alpha)).$ Now since  $ \oplus_{\alpha}B(L^{0}(\Omega_\alpha,
H_\alpha))$ is *-isomorphic to
$B(\oplus_{\alpha}L^{0}(\Omega_\alpha, H_\alpha)),$ one has that
the algebra  $L(M, \tau)$ is *-isomorphic with the algebra
$B(\oplus_{\alpha}L^{0}(\Omega_\alpha, H_\alpha)).$ The proof is
complete. $\blacksquare$

\textbf{Lemma 3.8}. \emph{If  $M$ is a von Neumann algebra of type
I, then $\mbox{mix}(S_0(M, \tau))$ is a standard subalgebra in $L(M,
\tau).$ }

 Proof. Let   $x\in L(M),
\tau)$ be a finite-generated operator. Then $q_\alpha x$ is a
finite-generated operator in  $L(L^{\infty}(\Omega_\alpha,
B(H_\alpha), \tau)$ for all $\alpha.$ By Proposition 3.4 one has
$q_\alpha x\in S_0(L^{\infty}(\Omega_\alpha, B(H_\alpha), \tau)$
for all $\alpha.$ Therefore
$x=\sum\limits_{\alpha}q_{\alpha}x\in\mbox{mix}(S_0(M, \tau)).$
The proof is complete. $\blacksquare$

Proof of Theorem 3.6.

Let  $d:S_0(M, \tau)\rightarrow S_0(M, \tau)$ be a $Z$-linear
derivation. Similar to the proof of Proposition 3.5  $d$ can be
extended to a  $Z_0$-linear derivation $\tilde{d}$ on
$\mbox{mix}(S_0(M, \tau)).$ Since  $\mbox{mix}(S_0(M, \tau))$ is a
standard algebra Theorem  2.1 implies that  $\tilde{d},$ and hence
$d,$ is implemented by an element of $L(M,\tau).$

\textbf{Remark 3.9.} The main result and the note after
Proposition 3.1 show that a derivation on $S_0(M, \tau)$  is
spatial if and only if it is $Z$-linear. Moreover [1, Example 4.6]
gives an example of non $Z$-linear (and hence non spatial)
derivation on $S_0(M, \tau)=L(M, \tau)$ for an appropriate von
Neumann algebra $M$ with a faithful normal finite trace $\tau.$

On the other hand if the lattice of projections in a von Neumann
algebra $M$ is atomic then any derivation on $S_0(M, \tau)$ is
automatically $Z$-linear (cf. [1, Corollary 4.7]).

Therefore we have

\textbf{Corollary 3.10.} \emph{If  $M$ is a von Neumann algebra
with the atomic lattice of projections, then any  derivation on
the algebra $S_0(M, \tau)$ is spatial, and in particular it is
continuous in the measure topology.}

\vspace{1cm}

\textbf{Acknowledgments.} \emph{The second and third named authors
would like to acknowledge the hospitality of the $\,$ "Institut
f\"{u}r Angewandte Mathematik",$\,$ Universit\"{a}t Bonn
(Germany). This work is supported in part by the DFG 436 USB
113/10/0-1 project (Germany) and the Fundamental Research
Foundation of the Uzbekistan Academy of Sciences.}

\end{document}